\noindent\centerline{\bf Erd\'elyi-Kober Fractional Integral Operators from a Statistical Perspective -II}

\vskip.3cm\centerline{A.M. MATHAI}
\vskip.2cm\centerline{Centre for Mathematical Sciences,}
\vskip.1cm\centerline{Arunapuram P.O., Pala, Kerala-68674, India, and}
\vskip.1cm\centerline{Department of Mathematics and Statistics, McGill University,}
\vskip.1cm\centerline{Montreal, Quebec, Canada, H3A 2K6}
\vskip.2cm\centerline{and}
\vskip.2cm\centerline{H.J. HAUBOLD}
\vskip.1cm\centerline{Office for Outer Space Affairs, United Nations}
\vskip.1cm\centerline{P.O. Box 500, Vienna International Centre} 
\vskip.1cm\centerline{A - 1400 Vienna, Austria, and}
\vskip.2cm\centerline{Centre for Mathematical Sciences,}
\vskip.1cm\centerline{Arunapuram P.O., Pala, Kerala-68674, India}

\vskip.5cm\noindent{\bf Abstract}

\vskip.3cm In this article we examine the densities of a product and a ratio of two real positive definite matrix-variate random variables $X_1$ and $X_2$, which are statistically independently distributed, and we consider the density of the product $U_1=X_2^{1\over2}X_1X_2^{1\over2}$ as well as the density of the ratio $U_2=X_2^{1\over2}X_1^{-1}X_2^{1\over2}$. We define matrix-variate Kober fractional integral operators of the first and second kinds from a statistical perspective, making use of the derivation in the predecessor of this paper for the scalar variable case, by deriving the densities of products and ratios where one variable has a matrix-variate type-1 beta density and the other variable has an arbitrary density. Various types of generalizations are considered, by using pathway models, by appending matrix variate hypergeometric series etc. During this process matrix-variate Saigo operator and other operators are also defined and properties studied.

\vskip.3cm\noindent{\bf 1.\hskip.3cm Introduction}

\vskip.3cm All the matrices appearing in this article are $p\times p$ real positive definite. Corresponding results for hermitian cases can be obtained but not given in this article. The following standard notations will be used. $A>O$ means the $p\times p$ real matrix is symmetric, $X=X'$, and further, it is positive definite. $|A|$ means the determinant of $A$, ${\rm tr}(A)= $ trace of $A$, ${\rm d}X$ will stand for the wedge product of differentials in any matrix $X$. If $X$ is $p\times q$ with $X=(x_{ij})$ then

$$\eqalignno{{\rm d}X&={\rm d}x_{11}\wedge{\rm d}x_{12}\wedge...\wedge{\rm d}x_{pq}\hbox{ for a general matrix}&(1.1)\cr
&=\prod_{i\ge j}\wedge{\rm d}x_{ij}=\prod_{j\ge i}\wedge{\rm d}x_{ij}\hbox{  when $X=X'$ or when $X$ is symmetric}.&(1.2)\cr}$$Also $\int_Xf(X){\rm d}X$ will mean the integrals over all $X$ (need not be symmetric or even square) of a real-valued scalar function $f(X)$ of matrix argument $X$. In the same format
$$\int_A^Bf(Y){\rm d}Y=\int_{O<A<Y<B}f(Y){\rm d}Y
$$will mean the integration of a real-valued scalar function of the real positive definite $p\times p$ matrix $Y$ over the space of positive definite matrices such that $A>O,X>O, X-A>O,B-X>O$ where $A$ and $B$ are $p\times p$ constant matrices. The notation will then imply that if $O<X<I$ then all eigenvalues of $X$ are in the open interval $(0,1)$. We will need some Jacobians of matrix transformations in this paper. For further results on Jacobians, see Mathai (1997).

$$\eqalignno{Y=AXA',X=X',|A|\ne 0&\Rightarrow {\rm d}Y=|A|^{p+1}{\rm d}X&(1.3)\cr
Y=X^{-1}&\Rightarrow {\rm d}Y=\cases{|X|^{-2p}{\rm d}X\hbox{ for a general $X$}\cr
|X|^{-(p+1)}{\rm d}X\hbox{ for $X=X'$}.\cr}&(1.4)\cr}
$$We will denote the unique positive definite square root of a positive definite matrix $A$ by $A^{1\over2}$. The following standard property will be used very often in this article. For $p\times p$ nonsingular matrices $A$ and $B$

$$\eqalignno{|I\pm AB|&=|I\pm BA|=|A|~|A^{-1}\pm B|=|B|~|B^{-1}\pm A|\hbox{  (when nonsingular)}\cr
|I\pm AB|&=|I\pm A^{1\over2}BA^{1\over2}|=|I\pm B^{1\over2}AB^{1\over2}|\hbox{ (when positive definite)}.&(1.5)\cr}
$$The real matrix-variate gamma function, denoted by $\Gamma_p(\alpha)$ is defined as follows which has an integral representation when $\Re(\alpha)>{{p-1}\over2}$.

$$\eqalignno{\Gamma_p(\alpha)&=\pi^{{p(p-1)}\over4}\Gamma(\alpha)\Gamma(\alpha-{1\over2})...\Gamma(\alpha-{{p-1}\over2}),\Re(\alpha)>{{p-1}\over2}\cr
&=\int_{X>O}|X|^{\alpha-{{p+1}\over2}}{\rm e}^{-{\rm tr}(X)}{\rm d}X,~\Re(\alpha)>{{p-1}\over2}.&(1.6)\cr}
$$The real matrix-variate type-1 beta density for the $p\times p$ positive definite matrix $X_1$, with parameters $\alpha$ and $\beta$ and denoted by $f_1(X_1)$, is defined as follows:

$$\eqalignno{f_1(X_1)&={{\Gamma_p(\alpha+\beta)}\over{\Gamma_p(\alpha)\Gamma_p(\beta)}}|X_1|^{\alpha-{{p+1}\over2}}
|I-X_1|^{\beta-{{p+1}\over2}},~O<X_1<I\cr
f_2(Y_1)&={{\Gamma_p(\alpha+\beta)}\over{\Gamma_p(\alpha)\Gamma_p(\beta)}}|Y_1|^{\beta-{{p+1}\over2}}
|I-Y_1|^{\alpha-{{p+1}\over2}}{\rm d}Y_1,O<Y_1<I\cr}
$$for $\Re(\alpha)>{{p-1}\over2},~\Re(\beta)>{{p-1}\over2}$, and $f(X_1)=0$, $f_2(Y_1)=0$ elsewhere. Type-1 and Type-2 beta integrals and beta functions are defined and denoted as follows for $\Re(\alpha)>{{p-1}\over2},\Re(\beta)>{{p-1}\over2}$:
$$\eqalignno{B_p(\alpha,\beta)&={{\Gamma_p(\alpha)\Gamma_p(\beta)}\over{\Gamma_p(\alpha+\beta)}}\cr
&=\int_{O<X<I}|X|^{\alpha-{{p+1}\over2}}|I-X|^{\beta-{{p+1}\over2}}{\rm d}X\hbox{ (type-1)}\cr
&=\int_{O<Y<I}|Y|^{\beta-{{p+1}\over2}}|I-Y|^{\alpha-{{p+1}\over2}}{\rm d}Y\hbox{ (type-1)}\cr
&=\int_{U>O}|U|^{\alpha-{{p+1}\over2}}|I+U|^{-(\alpha+\beta)}{\rm d}U\hbox{ (type-2)}\cr
&=\int_{V>O}|V|^{\beta-{{p+1}\over2}}|I+V|^{-(\alpha+\beta)}{\rm d}V\hbox{ (type-2)}&(1.7)\cr}
$$

\vskip.3cm\noindent{\bf 2.\hskip.3cm Kober Operator of the Second Kind for the Real Matrix-variate Case}

\vskip.3cm\noindent{\bf Definition 2.1.}\hskip.3cm We will define and denote Kober operator of the second kind for the real matrix-variate case as follows:

$$K_{X}^{\zeta,\alpha}f(X)={{|X|^{\zeta}}\over{\Gamma_p(\alpha)}}\int_{T>X}|T-X|^{\alpha-{{p+1}\over2}}
|T|^{-\zeta-\alpha}f(T){\rm d}T, \Re(\alpha)>{{p-1}\over2}.\eqno(2.1)
$$Consider two $p\times p$ real matrix-variate random variables $X_1$ and $X_2$, independently distributed, where $X_1$ has a real matrix-variate type-1 beta density $f_1(X_1)$ with parameters $(\zeta+{{p+1}\over2},\alpha)$, that is,

$$f_1(X_1)={{\Gamma_p(\zeta+\alpha+{{p+1}\over2})}\over{\Gamma_p(\zeta+{{p+1}\over2})\Gamma_p(\alpha)}}
|X_1|^{\zeta}|I-X_1|^{\alpha-{{p+1}\over2}},O<X_1<I
$$for $\Re(\zeta)>-1,\Re(\alpha)>{{p-1}\over2}$ and $f_1(X_1)=0$ elsewhere. Let $X_2$ have an arbitrary density $f(X_2)$. Then the joint density of $X_1$ and $X_2$ is $f_1(X_1)f(X_2)$. Let us consider the transformation $U=X_2^{1\over2}X_1X_2^{1\over2},X_2=V$ so that $X_2=V,X_1=V^{-{1\over2}}UV^{-{1\over2}}$ and the Jacobian is given by
${\rm d}X_1\wedge{\rm d}X_2=|V|^{-{{p+1}\over2}}{\rm d}U\wedge{\rm d}V$. If the joint density is denoted by $f(U,V)$ then
$$\eqalignno{f(U,V){\rm d}U\wedge{\rm d}V&={{\Gamma_p(\zeta+\alpha+{{p+1}\over2})}\over{\Gamma_p(\alpha)\Gamma(\zeta+{{p+1}\over2})}}
|V^{-{1\over2}}UV^{-{1\over2}}|^{\zeta}\cr
&\times|I-V^{-{1\over2}}UV^{-{1\over2}}|^{\alpha-{{p+1}\over2}}f(V)|V|^{-{{p+1}\over2}}{\rm d}U\wedge{\rm d}V.\cr}
$$Therefore the marginal density of $U$, denoted by $g(U)$, is available by integrating out $V$ from $f(U,V)$. That is,

$$\eqalignno{g(U)&=\int_Vf_1(V^{-{1\over2}}UV^{-{1\over2}})f(V)|V|^{-{1\over2}}{\rm d}V\cr
&={{\Gamma_p(\zeta+\alpha+{{p+1}\over2})}\over{\Gamma_p(\zeta+{{p+1}\over2})}}
\int_{V>U}{{1}\over{\Gamma_p(\alpha)}}|V|^{-{{p+1}\over2}}\cr
&\times |U|^{\zeta}|V|^{-\zeta}|V|^{-\alpha+{{p+1}\over2}}|V-U|^{\alpha-{{p+1}\over2}}f(V){\rm d}V\cr
&={{\Gamma_p(\alpha+\zeta+{{p+1}\over2})}\over{\Gamma_p(\zeta+{{p+1}\over2})}}K_{U}^{\zeta,\alpha}f(U).\cr}
$$Hence we have the following theorem:

\vskip.3cm\noindent{\bf Theorem 2.1.}\hskip.3cm{\it When $X_1$ and $X_2$ are independently distributed $p\times p$ positive definite real matrix random variables and when $X_2=V$ and $U=X_2^{1\over2}X_1X_2^{1\over2}$ or $X_1=V^{-{1\over2}}UV^{-{1\over2}}$ and when $X_1$ has a real matrix-variate type-1 beta distribution with the parameters $(\zeta+{{p+1}\over2},\alpha)$ and if $g(U)$ denotes the density of $U$ then

$${{\Gamma_p(\zeta+{{p+1}\over2})}\over{\Gamma_p(\alpha+\zeta+{{p+1}\over2})}}g(U)=K_{U}^{\zeta,\alpha}f(U)\eqno(2.2)$$is Kober fractional integral operator of the second kind for the real matrix-variate case.}

\vskip.3cm As a special case of (2.2), or independently, we can derive a result for the right sided Weyl operator for the real matrix-variate case. Let the right sided Weyl fractional integral operator be denoted by $({_XW_{\infty}^{-\alpha}}f)(X)$ with infinity here signifying that $T-X$ is positive definite.

\vskip.3cm\noindent{\bf Theorem 2.2.}\hskip.3cm{\it Let $X_1,X_2,U,V$ be as defined in Theorem 2.1. Let $X_1$ have a type-1 beta density with the parameters $({{p+1}\over2},\alpha)$. Let the arbitrary density of $X_2$ be denoted by $f_2(X_2)=|X_2|^{\alpha}f(X_2)$, where $f(X_2)$ is arbitrary. Let the density of $U$ be again denoted by $g(U)$. Then

$$\eqalignno{({_XW_{\infty}^{-\alpha}}f)(X)&={{1}\over{\Gamma(\alpha)}}\int_{V>U}|V-U|^{\alpha-{{p+1}\over2}}f(V){\rm  d}V\cr
&={{\Gamma_p({{p+1}\over2})}\over{\Gamma_p(\alpha+{{p+1}\over2})}}g(U),~ \Re(\alpha)>{{p-1}\over2}.&(2.3)\cr}
$$}

\vskip.3cm\noindent{\bf 2.1.\hskip.3cm A pathway generalization of Kober operator of the second kind in the matrix case}
\vskip.3cm A pathway generalization, parallel to the results in the scalar case can be considered. In the pathway case when generalization to matrices is considered we take $\delta=1$. For a general $\delta$, there will be problems with Jacobians of transformations for $X^{\delta}$ even if $X>O$ and $\delta>0$, see for example Mathai (1997) for the case $\delta=2$ and when $X=X'$ to see the type of complications. Hence we take the case $\delta=1$ only. Let $X_1$ have a pathway density
$$f_1(X_1)=C_1|X_1|^{\gamma}|I-a(1-q)X|^{{\eta}\over{1-q}}\eqno(2.4)
$$for $I-a(1-q)X>O, q<1,\eta>0,a>0$ where $C_1$ can be seen to be the following:
$$C_1={{[a(1-q)]^{p\gamma+{{p(p+1)}\over2}}\Gamma_p(\gamma+{{\eta}\over{1-q}}+(p+1))}
\over{\Gamma_p(\gamma+{{p+1}\over2})\Gamma_p({{\eta}\over{1-q}}+{{p+1}\over2})}}.\eqno(2.5)
$$Let $X_2$ have an arbitrary density $f(X_2)$ and let $X_1$ and $X_2$ be statistically independently distributed. Let $U=X_2^{1\over2}X_1X_2^{1\over2},X_2=V$ or $X_1=V^{-{1\over2}}UV^{-{1\over2}}$ and the Jacobian is $|V|^{-{{p+1}\over2}}$. Let $g(U)$ be the density of $U$. Then, going through the earlier steps we have the following:

$$g(U)=C_1|U|^{\gamma}\int_{V>a(1-q)U}|V|^{-\gamma-({{\eta}\over{1-q}}+{{p+1}\over2})}|V-a(1-q)U|^{{\eta}\over{1-q}}f(V){\rm d}V.\eqno(2.6)
$$Then
$$\eqalignno{\Gamma_p(\gamma+{{p+1}\over2})g(U)&={{[a(1-q)]^{p\gamma+{{p(p+1)}\over2}}
\Gamma_p(\gamma+{{\eta}\over{1-q}}+{{(p+1)}\over2})}\over{\Gamma_p({{\eta}\over{1-q}}+{{p+1}\over2})}}|U|^{\gamma}\cr
&\times\int_{V>a(1-q)U}|V-a(1-q)U|^{{\eta}\over{1-q}}|V|^{-\gamma-({{\eta}\over{1-q}}+{{p+1}\over2})}f(V){\rm d}V\cr
&=K_{U,a,q}^{\gamma,{{\eta}\over{1-q}}+{{p+1}\over2}}f(U)&(2.7)\cr}
$$where $K_{U,a,q}^{\gamma,{{\eta}\over{1-q}}+{{p+1}\over2}}f(U)$ can be called the generalized pathway Kober operator of the second kind in the real matrix-variate case. When the pathway parameter $q$ varies from $-\infty$ to 1 it provides a pathway or a class of operators and all these operators in this pathway class will eventually go to the exponential form. For $a=1,q=0,{{\eta}\over{1-q}}=\alpha -{{p+1}\over2}$ and $\gamma=\zeta$ we have

$$K_{U,a,q}^{\gamma,{{\eta}\over{1-q}}+{{p+1}\over2}}f(U)=K_{U}^{\zeta,\alpha}f(U)\eqno(2.8)
$$the Kober operator of the second kind as a constant multiple of the density of the product of two matrix-variate independently distributed random variables. Note that when $q\to 1_{-}$ we can evaluate the limit of $g(U)$ by using the following lemmas:

\vskip.3cm\noindent{\bf Lemma 2.1.}\hskip.3cm {\it
$$\lim_{q\to 1_{-}}C_1={{(a\eta)^{p\gamma+{{p(p+1)}\over2}}}\over{\Gamma_p(\gamma+{{p+1}\over2})}}\eqno(i)
$$}
\vskip.3cm\noindent{\bf Proof:}\hskip.3cm Open up each $\Gamma_p(\cdot)$ in $C_1$ of (2.5) in terms of ordinary gamma functions. Then use the following asymptotic approximation for gamma functions. For $|z|\to\infty$ and $\gamma$ a bounded quantity
$$\Gamma(z+\gamma)\approx \sqrt{2\pi}z^{z+\gamma-{1\over2}}{\rm e}^{-z}.\eqno(ii)
$$This is the first term in the asymptotic series. This term is also known as Stirling's approximation. When $q\to 1_{-}$ we have ${{1}\over{1-q}}\to \infty$ and hence take $|z|$ as ${{\eta}\over{1-q}}$ and expand all gammas by using Stirling's approximation to see that $C_1$ reduces to $(i)$ above.

\vskip.3cm\noindent{\bf Lemma 2.2.}\hskip.3cm{\it
$$\lim_{q\to 1_{-}}|I-a(1-q)X|^{{\eta}\over{1-q}}={\rm e}^{-a\eta~{\rm tr}(X)}.\eqno(iii)
$$}
\vskip.3cm\noindent{\bf Proof:}\hskip.3cm Writing the determinant in terms of eigenvalues we have
$$\eqalignno{|I-a(1-q)X|^{{\eta}\over{1-q}}&=\prod_{j=1}^p(1-a(1-q)\lambda_j)^{{\eta}\over{1-q}}&(iv)\cr
\noalign{\hbox{where $\lambda_1,...,\lambda_p$ are the eigenvalues of $X$. Now}}
\lim_{q\to 1_{-}}(1-a(1-q)\lambda_j)^{{\eta}\over{1-q}}&={\rm e}^{-a\eta~\lambda_j}.&(v)\cr
\noalign{\hbox{Hence}}
\lim_{q\to 1_{-}}|I-a(1-q)X|^{{\eta}\over{1-q}}&=\prod_{j=1}^p{\rm e}^{-a\eta~\lambda_j}={\rm e}^{-a\eta~{\rm tr}(X)}\cr}
$$which establishes (iii). Now by using Lemmas 2.1 and 2.2 we have

$$\eqalignno{\lim_{q\to 1_{-}}g(U)&={{(a\eta)^{p\gamma+{{p(p+1)}\over2}}}\over{\Gamma_p(\gamma+{{p+1}\over2})}}|U|^{\gamma}\cr
&\times\int_{V>O}|V|^{-\gamma-{{p+1}\over2}}{\rm e}^{-a\eta~{\rm tr}(V^{-{1\over2}}UV^{-{1\over2}})}{\rm d}V.&(2.9)\cr}
$$This is the limiting form of the pathway Kober operator of the second kind in this class of pathway operators of the second kind.
\vskip.2cm In the pathway generalization, one can also replace the parameter $a$ with a constant positive definite matrix $A$. In this case the model will be written as

$$f_1(X_1)=C_1(A)|X_1|^{\gamma}|I-(1-q)A^{1\over2}X_1A^{1\over2}|^{{\eta}\over{1-q}}\eqno(2.10)
$$for $q<1,A>O,X_1>O,I-(1-q)A^{1\over2}X_1A^{1\over2}>O$. The pathway parameter is still $q$. In this case

$$C_1(A)={{(1-q)^{p\gamma+{{p(p+1)}\over2}}|A|^{\gamma+{{p+1}\over2}}\Gamma_p(\gamma+{{\eta}\over{1-q}}+(p+1))}
\over{\Gamma_p(\gamma+{{p+1}\over2})\Gamma_p({{\eta}\over{1-q}}+{{p+1}\over2})}}.\eqno(2.11)
$$Then $g(U)$ of (2.6) goes to the following form, denoted by $g_A(U)$

$$\eqalignno{g_A(U)&=C_1(A)|A|^{{\eta}\over{1-q}}|U|^{\gamma}\cr
&\times\int_{V^{*}>O}|V^{1\over2}A^{-1}V^{1\over2}-(1-q)U|^{{\eta}\over{1-q}}|V|^{-\gamma-({{\eta}\over{1-q}}+{{p+1}\over2})}f(V){\rm d}V&(2.12)\cr
\noalign{\hbox{where}}
V^{*}&=V^{1\over2}A^{-1}V^{1\over2}-(1-q)U.\cr}
$$Then one can define a pathway generalized Kober operator of the second kind as

$$\eqalignno{K_{U,A,q}^{\gamma,{{\eta}\over{1-q}}+{{p+1}\over2}}f(U)&=\Gamma_p(\gamma+{{p+1}\over2})g_A(U)\cr
&={{(1-q)^{p\gamma+{{p(p+1)}\over2}}|A|^{\gamma+{{\eta}\over{1-q}}+{{p+1}\over2}}\Gamma_p(\gamma+{{\eta}\over{1-q}}
+(p+1))}\over{\Gamma_p({{\eta}\over{1-q}}+{{p+1}\over2})}}|U|^{\gamma}\cr
&\times\int_{V^{*}>O}|V^{1\over2}A^{-1}V^{1\over2}-(1-q)U|^{{\eta}\over{1-q}}
|V|^{-\gamma-({{\eta}\over{1-q}}+{{p+1}\over2})}f(V){\rm d}V&(2.13)\cr
\noalign{\hbox{In this case, as $q\to 1_{-}$ we have}}
\lim_{q\to 1_{-}}g_A(U)&={{|A|^{\gamma+{{p+1}\over2}}\eta^{p\gamma+{{p(p+1)}\over2}}}\over{\Gamma_p(\gamma+{{p+1}\over2})}}
|U|^{\gamma}\cr
&\times \int_{V>A}|V|^{-\gamma-{{p+1}\over2}}{\rm e}^{-\eta~{\rm tr}(A^{1\over2}V^{-{1\over2}}UV^{-{1\over2}}A^{1\over2})}f(V){\rm d}V.&(2.14)\cr}
$$
\vskip.3cm\noindent{3.\hskip.3cm M-transforms of Kober Operator of the Second Kind}
\vskip.3cm The generalized matrix transform or M-transform is defined and illustrated in Mathai (1997). The M-transform of Kober operator of the second kind is the following:

\vskip.3cm\noindent{\bf Theorem 3.1.}\hskip.3cm{\it For the Kober operator of the second kind defined in (2.1) the M-transform with parameter $s$ is given by
$$\eqalignno{M\{K_{X}^{\zeta,\alpha}f(X);s\}&=\int_{X>O}|X|^{s-{{p+1}\over2}}[\int_{T>X}{{|X|^{\zeta}}\over{\Gamma_p(\alpha)}}
|T-X|^{\alpha-{{p+1}\over2}}|T|^{-\zeta-\alpha}f(T){\rm d}T]{\rm d}X\cr
&={{\Gamma_p(\zeta+s)}\over{\Gamma(\alpha+\zeta+s)}}f^{*}(s),~\Re(\zeta+s)>{{p-1}\over2},\Re(\alpha)>{{p-1}\over2}&(3.1)\cr}
$$where $f^{*}(s)$ is the M-transform of $f(X)$.}

\vskip.3cm\noindent{\bf Proof:}\hskip.3cm Interchanging the integral we have
$$M\{K_{X}^{\zeta,\alpha}f(X);s\}=\int_{T>O}|T|^{-\zeta-\alpha}f(T)[{{1}\over{\Gamma_p(\alpha)}}
\int_{O<X<T}|X|^{\zeta+s-{{p+1}\over2}}|T-X|^{\alpha-{{p+1}\over2}}{\rm d}X]{\rm d}T.
$$Note that
$$|T-X|=|T|~|I-T^{-{1\over2}}XT^{-{1\over2}}|,~Y=T^{-{1\over2}}XT^{-{1\over2}}\Rightarrow{\rm d}Y=|T|^{-{{p+1}\over2}}{\rm d}X.
$$Hence
$$\int_{X<T}|X|^{\zeta-{{p+1}\over2}}|T-X|^{\alpha-{{p+1}\over2}}{\rm d}X=|T|^{\alpha+\zeta+s-{{p+1}\over2}}\int_Y|Y|^{\zeta+s-{{p+1}\over2}}|I-Y|^{\alpha-{{p+1}\over2}}{\rm d}Y.
$$We can evaluate the $Y$-integral by using real matrix-variate type-1 beta integral.
$$\int_{O<Y<I}|Y|^{\zeta+s-{{p+1}\over2}}|I-Y|^{\alpha-{{p+1}\over2}}{\rm d}Y={{\Gamma_p(\zeta+s)\Gamma_p(\alpha)}\over{\Gamma_p(\alpha+\zeta+s)}}
$$for $\Re(\alpha)>{{p-1}\over2},\Re(\zeta+s)>{{p-1}\over2}$. Now the $T$-integral gives
$$\int_{T>O}|T|^{s-{{p+1}\over2}}f(T){\rm d}T=f^{*}(s)
$$where $f^{*}(s)$ is the M-transform of $f(X)$. Hence (3.1) follows. Note that for $p=1$ the result agrees with that in the scalar case, which is available in the literature, see for example Mathai and Haubold (2008).
\vskip.2cm From (3.1) for $\zeta=0$ and $\Re(\alpha)>{{p-1}\over2}$ we have the special case of the Kober operator of the second kind
$$K_{X}^{0,\alpha}f(X)={{1}\over{\Gamma(\alpha)}}\int_{T>X}|T-X|^{\alpha-{{p+1}\over2}}|T|^{-\alpha}f(T){\rm d}T.\eqno(3.2)
$$But the right side of (3.2) is Weyl fractional integral of order $\alpha$ in the matrix case, ${_XW_{\infty}^{-\alpha}}f(X)$, except for the factor $|T|^{-\alpha}$. The Weyl integral in the matrix case is
$${_XW_{\infty}^{-\alpha}}f(X)={{1}\over{\Gamma_p(\alpha)}}\int_{T>X}|T-X|^{\alpha-{{p+1}\over2}}f(T){\rm d}T,\Re(\alpha)>{{p-1}\over2}.\eqno(3.3)
$$Hence we have the following corollary.

\vskip.3cm\noindent{\bf Corollary 3.1.1.}\hskip.3cm{\it The M-transform of the right sided Weyl operator in the real matrix case is given by
$$M\{{_XW_{\infty}^{-\alpha}}|X|^{-\alpha}f(X);s\}={{\Gamma_p(s)}\over{\Gamma_p(\alpha+s)}}f^{*}(s)\eqno(3.4)
$$for $\Re(s)>{{p-1}\over2},\Re(\alpha)>{{p-1}\over2}$ where $f^{*}(s)$ is the M-transform of $f(X)$.}

\vskip.3cm\noindent The proof is parallel to that in Theorem 3.1. Let us see whether a Mellin convolution type formula holds for Kober operator of the second kind in the matrix case. Let

$$g(U)=\int_V|V|^{-{{p+1}\over2}}f_1(V^{-{1\over2}}UV^{-{1\over2}})f_2(V){\rm d}V\eqno(3.5)
$$where $f_1(X_1)$ is a type-1 matrix-variate beta density with parameters $(\zeta+{{p+1}\over2},\alpha)$. That is,
$$f_1(X_1)={{\Gamma_p(\alpha+\zeta+{{p+1}\over2})}\over{\Gamma_p(\alpha)\Gamma_p(\zeta+{{p+1}\over2})}}|X_1|^{\zeta}
|I-X_1|^{\alpha-{{p+1}\over2}},O<X_1<I\eqno(3.6)
$$for $\Re(\alpha)>{{p-1}\over2},\Re(\zeta)>-1$ and $f_1(X_1)=0$ elsewhere. Substituting (3.6) in (3.5) we have

$$\eqalignno{{{\Gamma_p(\zeta+{{p+1}\over2})}\over{\Gamma_p(\alpha+\zeta+{{p+1}\over2})}}g(U)
&={{1}\over{\Gamma_p(\alpha)}}\int_V|V|^{-{{p+1}\over2}}|U|^{\zeta}|V|^{-\zeta}\cr
&\times|I-V^{-{1\over2}}UV^{-{1\over2}}|^{\alpha-{{p+1}\over2}}f(V){\rm d}V\cr
&={{|U|^{\zeta}}\over{\Gamma_p(\alpha)}}\int_V|V|^{-\zeta-\alpha}|V-U|^{\alpha-{{p+1}\over2}}f(V){\rm d}V\cr
&={{|U|^{\zeta}}\over{\Gamma_p(\alpha)}}\int_{V>U}|V-U|^{\alpha-{{p+1}\over2}}|V|^{-\zeta-\alpha}f(V){\rm d}V\cr
&=K_{U}^{\zeta,\alpha}f(U)&(3.7)\cr}
$$which is the Kober operator of the second kind. Hence we have the following theorem:

\vskip.3cm\noindent{\bf Theorem 3.2.}\hskip.3cm{\it Kober operator of the second kind with real matrix argument can also be represented as a Mellin convolution type formula

$$K_{X}^{\zeta,\alpha}f(X)=\int_V|V|^{-{{p+1}\over2}}f_1(V^{-{1\over2}}XV^{-{1\over2}})f_2(V){\rm d}V
$$where $f_1(X_1)$ is a type-1 beta density with parameters $(\zeta+{{p+1}\over2},\alpha)$ and $f_2(V)$ is an arbitrary function or arbitrary density if the Kober operator is to be taken as a constant multiple of a statistical density.}

\vskip.3cm\noindent{4.\hskip.3cm Generalization in Terms of Hypergeometric Series for Kober Operator of the Second Kind in the Real Matrix Case}

\vskip.3cm For introducing hypergeometric series of matrix argument we will need the definitions, notation and lemmas. Hypergeometric functions of matrix argument are defined in terms of matrix-variate Laplace transforms, M-transforms and zonal polynomials. Explicit series form for all cases is available through the definition in terms of zonal polynomials and hence we will define in terms of zonal polynomials.

$$\eqalignno{{_rF_s}(Z)&={_rF_s}(a_1,...,a_r;b_1,...,b_s;Z)\cr
&=\sum_{k=0}^{\infty}\sum_K{{(a)_K...(a_r)_K}\over{(b_1)_K...(b_s)_K}}{{C_K(Z)}\over{k!}}&(4.1)\cr
\noalign{\hbox{where $K=(k_1,...,k_p),k_1+...+k_p=k$ is a partition of $k=0,1,2,...$}}
(a)_K&=\prod_{j=1}^p(a-{{j-1}\over2})_{k_j}, (b)_{k_j}=b(b+1)...(b+k_j-1),(b)_0=1,b\ne 0&(4.2)\cr}
$$and $C_K(Z)$ is a zonal polynomial of order $k$ and $Z$ is a $p\times p$ matrix. The series is defined for the real and complex matrices. Zonal polynomials are certain symmetric functions of the eigenvalues of $Z$. In our discussions, $Z$ will be real and positive definite. For more details about zonal polynomials see Mathai, Provost and Hayakawa (1995). The following basic results are needed in our discussions. A standard notation in this area is

$$\Gamma_p(\alpha,K)=\Gamma_p(\alpha)(\alpha)_K.\eqno(4.3)
$$The following basic results are needed in our discussion.

\vskip.3cm\noindent{\bf Lemma 4.1.}\hskip.3cm{\it
$$\int_O^I|X|^{\alpha-{{p+1}\over2}}|I-X|^{\beta-{{p+1}\over2}}C_K(TX){\rm d}X={{\Gamma_p(\alpha,K)\Gamma_p(\beta)}\over{\Gamma_p(\alpha+\beta,K)}}C_K(T)\eqno(4.4)
$$for $\Re(\alpha)>{{p-1}\over2},\Re(\beta)>{{p-1}\over2}$.}

\vskip.3cm\noindent{\bf Lemma 4.2.}\hskip.3cm{\it For $\Re(\alpha)>{{p-1}\over2}$, $A>O$, $S>O$

$$\int_{O<S<A}|S|^{\alpha-{{p+1}\over2}}C_K(ZS){\rm d}S={{\Gamma_p(\alpha,K)\Gamma_p({{p+1}\over2})}\over{\Gamma_p(\alpha+{{p+1}\over2},K)}}|A|^{\alpha}C_K(ZA).\eqno(4.5)
$$}
\vskip.3cm Let us assume that all the parameters $a_1,...,a_r,b_1,...,b_s$ are real and positive and let the argument matrices be $p\times p$ and positive definite. For $A>O$, let the density of $X_1$ be

$$\eqalignno{f_1(X_1)&={{1}\over{c_f}}{_rF_s}(a_1,...,a_r;b_1,...,b_s;AX_1)
|X_1|^{\zeta}|I-X_1|^{\alpha-{{p+1}\over2}}\cr
&={{1}\over{c_f}}\sum_{k=0}^{\infty}\sum_K{{(a_1)_K...(a_r)_K}\over{(b_1)_K...(b_s)_K}}{{1}\over{k!}}
 C_K(AX_1)|X_1|^{\zeta}|I-X_1|^{\alpha-{{p+1}\over2}}&(4.6)\cr}
$$where the normalizing constant $c_f$ is available by integrating out term by term with the help of Lemma 4.1. It will be available in terms of a ${_{r+1}F_{s+1}}$. Let $f(X_2)$ be an arbitrary density. As before, let $X_1=V^{-{1\over2}}UV^{-{1\over2}}$, then denoting the density of $U$ again by $g(U)$ we have

$$\eqalignno{g(U)&=\int_vf_1(V^{-{1\over2}}UV^{-{1\over2}})f(V)|V|^{-{{p+1}\over2}}{\rm d}V\cr
&={{1}\over{c_f}}{{\Gamma_p(\alpha+\zeta+{{p+1}\over2})}\over{\Gamma_p(\zeta+{{p+1}\over2})\Gamma_p(\alpha)}}(\sum_{k=0}^{\infty}\sum_K{{(a_1)_K...(a_r)_K}\over{(b_1)_K...(b_s)_K}}{{1}\over{k!}}\cr
&\times \int_V|V^{-{1\over2}}UV^{-{1\over2}}|^{\zeta}|I-V^{-{1\over2}}UV^{-{1\over2}}|^{\alpha-{{p+1}\over2}}
|V|^{-{{p+1}\over2}}C_K(AV^{-{1\over2}}UV^{-{1\over2}})f(V){\rm d}V&(4.7)\cr}
$$This is the generalization of a constant times the Kober operator of the second kind in the matrix case. For ${_rF_s}={_2F_1}$ one has the matrix-variate generalization of a constant times the Saigo operator of the second kind in the real matrix-variate case.

\vskip.3cm\noindent{\bf 5.\hskip.3cm Kober Fractional Integral Operators of the First Kind in the Matrix Case}

\vskip.3cm\noindent{\bf Definition 5.1.}\hskip.3cm We will give the following definition and notation for Kober fractional integral operator of the first kind in the real matrix-variate case:
$$I_{X}^{\zeta,\alpha}f(X)={{|X|^{-\zeta-\alpha}}\over{\Gamma_p(\alpha)}}\int_{V<X}|X-V|^{\alpha-{{p+1}\over2}}
|V|^{\zeta}f(V){\rm d}V\eqno(5.1)$$for $\Re(\zeta)>-1,\Re(\alpha)>{{p-1}\over2}$.

\vskip.3cm\noindent{\bf Theorem 5.1.}\hskip.3cm{\it For $\Re(\alpha)>{{p-1}\over2}, \Re(\zeta)>-1$ the M-transform, with parameter $s$, of Kober operator of the first kind in the real matrix-variate case, is given by

$$\eqalignno{M\{I_{X}^{\zeta,\alpha}f(X);s\}&=\int_{X>O}|X|^{s-{{p+1}\over2}}[{{|X|^{-\zeta-\alpha}}\over{\Gamma_p(\alpha)}}
\int_{V<X}|X-V|^{\alpha-{{p+1}\over2}}|V|^{\zeta}f(V){\rm d}V]{\rm d}X\cr
&={{\Gamma_p(\zeta+{{p+1}\over2}-s)}\over{\Gamma_p(\alpha+\zeta+{{p+1}\over2}-s)}}f^{*}(s),\Re(s)<\Re(\zeta+1),\Re(\alpha)>{{p-1}\over2}&(5.2)\cr}
$$where $f^{*}(s)$ is the M-transform of $f(X)$.}

\vskip.3cm\noindent{\bf Proof:}\hskip.3cm Integrating out $X$ first we have the $X$-integral
$$\eqalignno{\int_{X>V}|X|^{s-\zeta-\alpha-{{p+1}\over2}}|X-V|^{\alpha-{{p+1}\over2}}{\rm d}X&=\int_{Y>O}|Y+V|^{s-\zeta-\alpha-{{p+1}\over2}}|Y|^{\alpha-{{p+1}\over2}}{\rm d}Y,Y=X-V\cr
&=|V|^{s-\zeta-\alpha-{{p+1}\over2}}\int_{Y>O}|I+V^{-{1\over2}}YV^{-{1\over2}}|^{s-\zeta-\alpha-{{p+1}\over2}}
|Y|^{\alpha-{{p+1}\over2}}{\rm d}Y.\cr}
$$Put $Z=V^{-{1\over2}}YV^{-{1\over2}}\Rightarrow {\rm d}Z=|V|^{-{{p+1}\over2}}{\rm d}Y$. Then the $X$-integral is

$$\eqalignno{|V|^{s-\zeta-{{p+1}\over2}}&\int_{Z>O}|Z|^{\alpha-{{p+1}\over2}}|I+Z|^{-({{p+1}\over2}+\alpha+\zeta-s)}{\rm d}Z\cr
&=|V|^{s-\zeta-{{p+1}\over2}}{{\Gamma_p(\alpha)\Gamma_p({{p+1}\over2}+\zeta-s)}\over{\Gamma_p({{p+1}\over2}+\alpha+\zeta-s)}}\cr}
$$for $\Re(\alpha)>{{p-1}\over2},\Re(\zeta-s)>-1$ by evaluating the integral by using a type-2 matrix-variate beta integral in the real case. Now, the $V$-integral becomes
$$\int_{V>O}|V|^{s-{{p+1}\over2}}f(V){\rm d}V=f^{*}(s).
$$Hence
$$M\{I_{X}^{\zeta,\alpha}f(X);s\}={{\Gamma_p({{p+1}\over2}+\zeta-s)}\over{\Gamma_p({{p+1}\over2}+\alpha+\zeta-s)}}f^{*}(s)\eqno(5.3)
$$for $\Re(\alpha)>{{p-1}\over2},\Re(\zeta-s)>-1$. Note that for $\zeta=0$,
$$I_{X}^{0,\alpha}f(X)=|X|^{-\alpha}{_0D_X^{-\alpha}}f(X)\eqno(5.4)
$$where ${_0D_X^{-\alpha}}$ is the left sided Riemann-Liouville fractional integral for the matrix-variate case. Note that for the scalar case, for $p=1$,
$$M\{I_{x}^{\zeta,\alpha}f(x);s\}={{\Gamma(1+\zeta-s)}\over{\Gamma(1+\alpha+\zeta-s)}}
\eqno(5.5)
$$for $\Re(\alpha)>0,\Re(\zeta-s)>-1$ agreeing with the corresponding Mellin transform in the scalar case.

\vskip.3cm\noindent{\bf Corollary 5.1.1}.\hskip.3cm{\it The M-transform of $|X|^{-\alpha}{_0D_X^{-\alpha}}f(X)$ is given by
$$M\{|X|^{-\alpha}{_0D_X^{-\alpha}}f(X);s\}={{\Gamma_p({{p+1}\over2}-s)}\over{\Gamma_p({{p+1}\over2}+\alpha-s)}}f^{*}(s)\eqno(5.6)
$$for $\Re(\alpha)>{{p-1}\over2}, \Re(s)<1$.}

\vskip.3cm\noindent The proof is parallel to that in Theorem 5.1.
\vskip.2cm Let us treat a Kober operator of the first kind as a statistical density. Let $X_2$ have an arbitrary real matrix-variate density $f(X_2)$ and $X_1$ have a real matrix-variate type-1 beta density with parameters $(\zeta,\alpha)$. That is,

$$f_1(X_1)={{\Gamma_p(\zeta+\alpha)}\over{\Gamma_p(\zeta)\Gamma_p(\alpha)}}|X_1|^{\zeta-{{p+1}\over2}}
|I-X_1|^{\alpha-{{p+1}\over2}},O<X_1<I\eqno(5.7)
$$for $\Re(\zeta)>{{p-1}\over2},\Re(\alpha)>{{p-1}\over2}$ and $f_1(X_1)=0$ elsewhere. Let $X_1$ and $X_2$ be statistically independently distributed. Consider the transformation $X_2=V,X_1=V^{1\over2}U^{-1}V^{1\over2}$. The Jacobian is given by
$${\rm d}X_1\wedge{\rm d}X_2=|V|^{{p+1}\over2}|U|^{-(p+1)}{\rm d}U\wedge{\rm d}V.
$$The marginal density of $U$,  denoted by $g_r(U)$ where $r$ designates that it is coming from a ratio, is given by
$$\eqalignno{g_r(U)&={{\Gamma_p(\zeta+\alpha)}\over{\Gamma_p(\zeta)\Gamma_p(\alpha)}}
\int_V|V^{1\over2}U^{-1}V^{1\over2}|^{\zeta-{{p+1}\over2}}\cr
&\times|I-V^{1\over2}U^{-1}V^{1\over2}|^{\alpha-{{p+1}\over2}}f(V)|V|^{{p+1}\over2}|U|^{-(p+1)}{\rm d}V\cr
&={{\Gamma_p(\zeta+\alpha)}\over{\Gamma_p(\zeta)\Gamma_p(\alpha)}}|U|^{-\zeta-\alpha}
\int_{V<U}|U-V|^{\alpha-{{p+1}\over2}}|V|^{\zeta}f(V){\rm d}V.\cr
\noalign{\hbox{Therefore}}
{{\Gamma_p(\zeta)}\over{\Gamma_p(\zeta+\alpha)}}g_r(U)&={{|U|^{-\zeta-\alpha}}\over{\Gamma_p(\alpha)}}
\int_{V<U}|V|^{\zeta}f(V){\rm d}V\cr
&=I_{X}^{\zeta,\alpha}f(X).&(5.8)\cr}
$$This is Kober operator of the first kind in the real matrix-variate case and it can be considered as a constant multiple of a real matrix-variate statistical density. \vskip.2cm One can also consider a pathway extension for the real matrix-variate Kober operator of the first kind.

\vskip.3cm\noindent{\bf 5.2.\hskip.3cm Pathway Extension of Kober Operator of the First Kind in the Matrix Case}

\vskip.3cm Consider the following pathway modified form of the density for $X_1$. That is,

$$f_1(X_1)=C_2|X_1|^{\gamma-{{p+1}\over2}}|I-a(1-q)X_1|^{{\eta}\over{1-q}},I-a(1-q)X_1>O\eqno(5.9)
$$for $q<1,a>0,\eta>0$ where
$$C_2={{[a(1-q)]^{p\gamma}\Gamma_p(\gamma+{{\eta}\over{1-q}}+{{p+1}\over2})}\over{\Gamma_p({{\eta}\over{1-q}}+{{p+1}\over2})
\Gamma_p(\gamma)}}.\eqno(5.10)
$$Consider the same type of transformation as before: $X_2=V,X_1=V^{1\over2}U^{-1}V^{1\over2}$. The marginal density of $U$, denoted by $g_p(U)$, is given by

$$\eqalignno{g_p(U)&=C_2\int_V|V^{1\over2}U^{-1}V^{1\over2}|^{\gamma-{{p+1}\over2}}
|I-a(1-q)V^{1\over2}U^{-1}V^{1\over2}|^{{\eta}\over{1-q}}\cr
&\times f(V)|V|^{{p+1}\over2}|U|^{-(p+1)}{\rm d}V&(5.11)\cr
&=C_2|U|^{-\gamma-({{\eta}\over{1-q}}+{{p+1}\over2})}\int_{U>a(1-q)V}|U-a(1-q)V|^{{\eta}\over{1-q}}|V|^{\gamma}f(V){\rm d}V.\cr
\noalign{\hbox{Then}}
\Gamma_p(\gamma)g_p(U)&={{[a(1-q)]^{p\gamma}\Gamma_p(\gamma+{{\eta}\over{1-q}}+{{p+1}\over2})}
\over{\Gamma_p({{\eta}\over{1-q}}+{{p+1}\over2})}}|U|^{-\gamma-({{\eta}\over{1-q}}+{{p+1}\over2})}\cr
&\times\int_{U>a(1-q)V}|U-a(1-q)V|^{{\eta}\over{1-q}}|V|^{\gamma}f(V){\rm d}V&(5.12)\cr}
$$The right side of (5.12) is the pathway extension of Kober operator of the first kind. The right side divided by $\Gamma_p(\gamma)$ is also a statistical density of a type of ratio of independently distributed matrix-variate random variables.
\vskip.2cm Note that for $a=1,q=0,{{\eta}\over{1-q}}+{{p+1}\over2}=\alpha$, (5.12) reduces to the special case (5.7) for $\gamma=\zeta$. Thus, (5.12) describes a vast family of operators which can all be considered as generalizations of the Kober operator of the first kind in the real matrix-variate case. The limiting form when $q\to 1_{-}$ is available from the structure in (5.11). Note that
$$\eqalignno{\lim_{q\to 1_{-}}&|I-a(1-q)V^{1\over2}U^{-1}V^{1\over2}|^{{\eta}\over{1-q}}={\rm e}^{-a\eta~{\rm tr}(V^{1\over2}U^{-1}V^{1\over2})}&(5.13)\cr
\noalign{\hbox{Hence}}
\lim_{q\to 1_{-}}g_p(U)&=(\lim_{q\to 1_{-}}C_2)\int_{V>O}|U|^{-\gamma-{{p+1}\over2}}|V|^{\gamma}\cr
&\times {\rm e}^{-a\eta~{\rm tr}(V^{1\over2}U^{-1}V^{1\over2})}f(V){\rm d}V&(5.14)\cr
\noalign{\hbox{where}}
\lim_{q\to 1_{-}}C_2&={{(a\eta)^{p\gamma}}\over{\Gamma_p(\gamma)}}.\cr
\noalign{\hbox{That is,}}
\lim_{q\to 1_{-}}g_p(U)&={{(a\eta)^{p\gamma}}\over{\Gamma_p(\gamma)}}|U|^{-\gamma-{{p+1}\over2}}\cr
&\times\int_{V>O}|V|^{\gamma}{\rm e}^{-a\eta~{\rm tr}(V^{1\over2}U^{-1}V^{1\over2})}f(V){\rm d}V.&(5.15)\cr}
$$In this case also one can replace the parameter $a$ in $f_1(X_1)$ by a constant positive definite matrix $A>O$. Then $f_1(X_1)$ denoted by $f_{1A}(X_1)$ can be written as

$$\eqalignno{f_{1A}(X_1)&=C_2(A)|X_1|^{\gamma-{{p+1}\over2}}|I-(1-q)A^{1\over2}X_1A^{1\over2}|^{{\eta}\over{1-q}}\cr
\noalign{\hbox{where}}
C_2(A)&={{(1-q)^{p\gamma}|A|^{\gamma}\Gamma_p(\gamma+{{\eta}\over{1-q}}+{{p+1}\over2})}\over{\Gamma_p(\gamma)
\Gamma_p({{\eta}\over{1-q}}+{{p+1}\over2})}}.\cr}
$$Then the density of $U=V^{1\over2}U^{-1}V^{1\over2}$, denoted by $g_A(U)$, is given by

$$\eqalignno{g_A(U)&=C_2(A)|U|^{-\gamma-({{\eta}\over{1-q}}+{{p+1}\over2})}\cr
&\times\int_{U>(1-q)V^{1\over2}AV^{1\over2}}|V|^{\gamma}|U-(1-q)V^{1\over2}AV^{1\over2}|^{{\eta}\over{1-q}}f(V){\rm d}V.&(5.16)\cr}
$$

\vskip.3cm\noindent{\bf Acknowledgment}

\vskip.3cm The authors would like to thank the Department of Science and Technology, Government of India, New Delhi, for the financial assistance for this work under project number SR/S4/MS:287/05.

\vskip.3cm\centerline{\bf References}

\vskip.3cm\noindent Mathai, A.M. (1997):\hskip.3cm{\it Jacobians of Matrix Transformations and Functions of Matrix Argument}, World Scientific Publishing, New York.
\vskip.3cm\noindent Mathai, A.M. and Haubold, H.J. (2008):\hskip.3cm {\it Special Functions for Applied Scientists}, Springer, New York.
\vskip.3cm\noindent Mathai, A.M., Provost, S.B., and Hayakawa, T. (1995):\hskip.3cm{\it Bilinear Forms and Zonal Polynomials}, Springer, New York.

\bye